\documentclass[a4paper]{amsproc}

\usepackage{amssymb}
\usepackage{pstricks}
\usepackage{pstcol}
\usepackage{graphicx}
\usepackage{amsmath}
\usepackage{amsfonts}
\usepackage{latexsym}

\theoremstyle{plain}

\theoremstyle{definition}
\newtheorem{exm}{Example}[section]

\numberwithin{equation}{section}
\renewcommand{\leq}{\leqslant}
\renewcommand{\geq}{\geqslant}

\title[Frullani integrals]{The integrals in Gradshteyn and Ryzhik. \\
Part 15: Frullani integrals}

\subjclass[2000]{Primary 33}


\author[M. Albano]{Matthew Albano}
\address{Department of Mathematics,
New Jersey Institute of Technology, Newark, NJ 07102}
\email{rocky23cancook@gmail.com}

\author[T. Amdeberhan]{Tewodros Amdeberhan}
\address{Department of Mathematics,
Tulane University, New Orleans, LA 70118}
\email{tamdeber@math.tulane.edu}

\author[E. Beyerstedt]{Erin Beyerstedt}
\address{Department of Mathematics,
Tulane University, New Orleans, LA 70118}
\email{ebeyerst@math.tulane.edu}

\author[V. Moll]{Victor H. Moll}
\address{Department of Mathematics,
Tulane University, New Orleans, LA 70118}
\email{vhm@math.tulane.edu}

\address{\hfill{\it Received 4 20
2010 revised ?? }\newline Departamento de Matem\'atica
\newline
Universidad T\'ecnica Federico Santa Mar\'{\i}a
\newline  Casilla 110-V,
\newline Valpara\'{\i}so, Chile}

\begin{document}

{\begin{flushleft}\baselineskip9pt\scriptsize {\bf SCIENTIA}\newline
Series A: {\it Mathematical Sciences}, Vol. ?? (2010), ??
\newline Universidad T\'ecnica Federico Santa Mar{\'\i}a
\newline Valpara{\'\i}so, Chile
\newline ISSN 0716-8446
\newline {\copyright\space Universidad T\'ecnica Federico Santa
Mar{\'\i}a\space 2010}
\end{flushleft}}

\vspace{10mm} \setcounter{page}{1} \thispagestyle{empty}

\begin{abstract}
The table of Gradshteyn and Ryzhik contains some  integrals that 
can be reduced to the Frullani type. We present a selection of them. 
\end{abstract}

\maketitle

\newcommand{\nn}{\nonumber}
\newcommand{\ba}{\begin{eqnarray}}
\newcommand{\ea}{\end{eqnarray}}
\newcommand{\ift}{\int_{0}^{\infty}}
\newcommand{\ione}{\int_{0}^{1}}
\newcommand{\ifft}{\int_{- \infty}^{\infty}}
\newcommand{\no}{\noindent}
\newcommand{\realpart}{\mathop{\rm Re}\nolimits}
\newcommand{\imagpart}{\mathop{\rm Im}\nolimits}

\newtheorem{Definition}{\bf Definition}[section]
\newtheorem{Thm}[Definition]{\bf Theorem} 
\newtheorem{Example}[Definition]{\bf Example} 
\newtheorem{Lem}[Definition]{\bf Lemma} 
\newtheorem{Note}[Definition]{\bf Note} 
\newtheorem{Cor}[Definition]{\bf Corollary} 
\newtheorem{Prop}[Definition]{\bf Proposition} 
\newtheorem{Problem}[Definition]{\bf Problem} 
\numberwithin{equation}{section}

\maketitle

\section{Introduction} \label{intro} 
\setcounter{equation}{0}

The table of integrals \cite{gr} contains many evaluations of the form
\begin{equation}
\ift \frac{f(ax)-f(bx)}{x} \, dx = \left[ f(0) - f(\infty) \right] \, 
\ln \left(\frac{b}{a} \right).
\label{fru-1}
\end{equation}
\noindent
Expressions of this type are called {\em Frullani integrals}. Conditions 
that guarantee the validity of this formula are given in \cite{arias1} 
and \cite{ostrowski1}. In particular, the 
continuity of $f'$ and the convergence of the integral are sufficient for 
(\ref{fru-1}) to hold.

\section{A list of examples} \label{sec-examples} 
\setcounter{equation}{0}

Many of the entries in \cite{gr} are simply particular cases of (\ref{fru-1}). 

\begin{exm}
\label{example-1}
The evaluation of $3.434.2$ in \cite{gr}:
\begin{equation}
\ift \frac{e^{-ax} - e^{-bx}}{x} \, dx = \ln \left( \frac{b}{a} \right)
\end{equation}
\noindent
corresponds to the function $f(x) = e^{-x}$. 
\end{exm}

\begin{exm}
The change of variables 
$t = e^{-x}$ in Example \ref{example-1} yields
\begin{equation}
\int_{0}^{1} \frac{t^{b-1}-t^{a-1}}{\ln t} \, dt = 
\ln \left( \frac{a}{b} \right).
\end{equation}
\noindent
This is $4.267.8$ in \cite{gr}.
\end{exm}

\begin{exm}
A generalization of the previous example appears as entry $3.476.1$ 
in \cite{gr}:
\begin{equation}
\int_{0}^{\infty} \left(e^{-vx^{p}} - e^{-ux^{p}} \right)  \, \frac{dx}{x}= 
\frac{1}{p} \ln \left( \frac{u}{v} \right). 
\end{equation}
\noindent 
This comes from Frullani's result with a simple additional scaling. 
\end{exm}

\begin{exm}
The choice
\begin{equation}
f(x) = \frac{e^{-qx}-e^{-px}}{x},
\end{equation}
\noindent
with $p, \, q >0$ satisfies $f(\infty) = 0$ and 
\begin{equation}
f(0) = \lim\limits_{x \to 0} \frac{e^{-qx}-e^{-px}}{x} = p-q.
\end{equation}
\noindent
Then Frullani's theorem yields
\begin{equation}
\ift \left( \frac{e^{-aqx}-e^{-apx}}{ax} - \frac{e^{-bqx}-e^{-bpx}}{bx} \right)
\, \frac{dx}{x} = (p-q) \ln \left( \frac{b}{a} \right),
\nonumber 
\end{equation}
\noindent
that can be written as 
\begin{equation}
\ift \left( \frac{e^{-aqx}-e^{-apx}}{a} - \frac{e^{-bqx}-e^{-bpx}}{b} \right)
\, \frac{dx}{x^{2}} = (p-q) \ln \left( \frac{b}{a} \right).
\nonumber 
\end{equation}
\noindent
This is entry $3.436$ in \cite{gr}. \\
\end{exm}

\begin{exm}
Now choose 
\begin{equation}
f(x) = \frac{x}{1-e^{-x}} \text{exp}(-c e^{x}). 
\end{equation}
\noindent
Then Frullani's theorem yields entry $3.329$ of \cite{gr}, 
in view of $f(0) = e^{-c}$ and 
$f(\infty) = 0$:
\begin{equation}
\ift \left( \frac{a \, \text{exp}(-ce^{ax})}{1-e^{-ax}} - 
     \frac{b \, \text{exp}(-ce^{bx})}{1-e^{-bx}}   \right) \, dx = 
e^{-c} \, \ln \left( \frac{b}{a} \right). 
\nonumber
\end{equation}
\end{exm}

\medskip

\begin{exm}
The next example uses 
\begin{equation}
f(x) = (x+c)^{-\mu},
\end{equation}
\noindent
with $c, \, \mu > 0$, to produce
\begin{equation}
\ift \frac{(ax+c)^{-\mu} - (bx+c)^{-\mu}}{x} \, dx = 
c^{-\mu} \ln \left( \frac{b}{a} \right). 
\end{equation}
\noindent
This is $3.232$ in \cite{gr}. 
\end{exm}

\begin{exm}
Entry $4.536.2$ in \cite{gr} is
\begin{equation}
\int_{0}^{\infty} \frac{\tan^{-1}(px) - \tan^{-1}(qx)}{x} \, dx = 
\frac{\pi}{2} \ln \left( \frac{p}{q} \right). 
\end{equation}
\noindent
This follows directly from (\ref{fru-1}) by choosing $f(x) = \tan^{-1}x$. 
\end{exm}

\begin{exm}
The function $f(x) = \ln(a + be^{-x})$ gives the evaluation of entry 
$4.319.3$ of \cite{gr}:
\begin{equation}
\int_{0}^{\infty} 
\frac{\ln(a+be^{-px}) - \ln(a + be^{-qx})}{x} \, dx = 
\ln \left( \frac{a}{a+b} \right) \, \ln \left( \frac{p}{q} \right).
\end{equation}
\end{exm}

\begin{exm}
The function $f(x) = ab \ln(1 +x)/x$ produces entry 
$4.297.7$ of \cite{gr}:
\begin{equation}
\int_{0}^{\infty} 
\frac{b \ln(1 + a x) - a \ln(1+ bx)}{x^{2}} \, dx = ab \ln \left( \frac{b}{a} 
\right). 
\end{equation}
\end{exm}

\begin{exm}
Entry $3.484$:
\begin{equation}
\int_{0}^{\infty} \left[ \left( 1 + \frac{a}{qx} \right)^{qx} - 
\left( 1 + \frac{a}{px} \right)^{px} \right] \, \frac{dx}{x}  = 
\left( e^{a} - 1 \right) \, \ln \left( \frac{q}{p} \right),
\end{equation}
\noindent
is obtained by choosing $f(x) = (1 + a/x)^{x}$ in (\ref{fru-1}). 
\end{exm}

\begin{exm}
The final example in this section corresponds to the function 
\begin{equation}
f(x) = \frac{a + be^{-x}}{ce^{x} + g + he^{-x}}
\end{equation}
\noindent
that produces entry $3.412.1$ of \cite{gr}:
\begin{equation}
\ift \left[ \frac{a + be^{-px}}{ce^{px} + g + he^{-px}} - 
 \frac{a + be^{-qx}}{ce^{qx} + g + he^{-qx}}  \right] \, \frac{dx}{x} = 
\frac{a+b}{c+g+h} \ln \left( \frac{q}{p} \right). 
\end{equation}
\end{exm}

\section{A separate source of examples} \label{rama-1} 
\setcounter{equation}{0}

The list presented in this section contains integrals of Frullani type 
that were found in volume 1 of Ramanujan's Notebooks \cite{berndtI}. 

\begin{exm}
\begin{equation}
\int_0^{\infty} \frac{\tan^{-1} ax - \tan^{-1} bx}{x} \, dx \; = \; 
\frac{\pi}{2} \ln \frac{a}{b}
\nonumber
\end{equation}
\end{exm}

\begin{exm}
\begin{equation}
\int_0^{\infty} \ln \frac{p + qe^{-ax}}{p + qe^{-bx}} \frac{dx}{x} 
\; = \; \ln \left( 1 + \frac{q}{p} \right) \ln \frac{b}{a}
\nonumber
\end{equation}
\end{exm}

\begin{exm}
\begin{equation}
\int_0^{\infty} \left[ \left(\frac{ax+p}{ax+q}\right)^n - \left(\frac{bx+p}{bx+q}\right)^n \right] \frac{dx}{x} \; = \; \left( 1 - \frac{p^n}{q^n} \right) 
\ln \frac{a}{b}
\nonumber
\end{equation}
\noindent
where $a,b,p,q$ are all positive.
\end{exm}

\begin{exm}
\begin{equation}
\int_0^{\infty} \frac{\cos ax - \cos bx}{x} \, dx \; = \; \ln \frac{b}{a}
\nonumber
\end{equation}
\end{exm}

\begin{exm}
\begin{equation}
\int_0^{\infty} \sin \left( \frac{(b-a)x}{2} \right) \, 
\sin \left(\frac{(b+a)x}{2} \right) \, 
\frac{dx}{x} 
\; = \; 
\int_0^{\infty} \frac{\cos ax  - \cos bx }{2x} \, dx  \; = \; \frac{1}{2} 
\ln \frac{b}{a}
\nonumber
\end{equation}
\end{exm}

\begin{exm}
\begin{equation}
\int_0^{\infty} \sin px \, \sin qx \, \frac{dx}{x} \; = 
\; \int_0^{\infty} \frac{\cos[(p-q)x] - \cos[(p+q)x]}{2x} \, dx \; = \; \frac{1}{2} \ln \frac{p+q}{p-q}
\nonumber
\end{equation}
\end{exm}

\begin{exm}
The evaluation of 
\begin{equation}
\int_0^{\infty} \ln \left(\frac{1 + 2n \cos ax + n^2}{1 + 2n \cos bx + n^2}\right) \frac{dx}{x} \; = \; \left\{
\begin{array}{ll}
\ln \frac{b}{a} \ln (1 + n)^2 &  n^2 < 1\\
\ln \frac{b}{a} \ln \left(1 + \frac{1}{n}\right)^2 &  n^2 > 1
\end{array} \right.
\nonumber
\end{equation}
\noindent
is more delicate and is given in detail in the next section.
\end{exm}

\begin{exm}
The value 
\begin{equation}
\int_{0}^{\infty} \frac{e^{-ax} \sin ax - e^{-bx} \sin bx}{x} \, dx \; = 
\; 0 
\nonumber
\end{equation}
\noindent
follows directly from (\ref{fru-1}) since, in this case $f(x) = 
e^{-x} \sin x$ satisfies  $f(\infty) = f(0) = 0$.
\end{exm}

\begin{exm}
\begin{equation}
\int_0^{\infty} \frac{e^{-ax} \cos ax - e^{-bx} \cos bx}{x} \, dx \; = \; 
\ln \frac{b}{a}.
\nonumber
\end{equation}
\end{exm}

\section{A more delicate example} \label{delicate} 
\setcounter{equation}{0}

Entry $4.324.2$ of \cite{gr} states that 
\begin{multline}
\int_{0}^{\infty} \left[ \ln(1 + 2 a \cos px + a^{2} )  - 
\ln(1 + 2 a \cos qx + a^{2} )  \right] \, \frac{dx}{x} = 
\\
\begin{cases}
2 \ln \left( \frac{q}{p} \right) \ln(1+a) & \quad -1 < a \leq 1 \\
2 \ln \left( \frac{q}{p} \right) \ln(1 + 1/a) & \quad a< -1 \text{ or } a \geq 1.
\end{cases}
\end{multline}
\noindent
This requires a different approach since the obvious candidate for a direct 
application of Frullani's theorem, namely $f(x) = \ln(  1 + 2a \cos x + 
a^{2} )$, does not have a limit at infinity. 

In order to evaluate this entry, start with
\begin{equation}
\int_0^1x^ydx=\frac1{y+1},
\end{equation}
\noindent
so
\begin{equation}
\int_0^1dy\int_0^1x^ydx=\int_0^1dx\int_0^1x^ydy=
\int_0^1\frac{x-1}{\ln x}dx=\int_0^1\frac{dy}{y+1}=\ln 2.
\end{equation}
This is now generalized for arbitrary symbols $\alpha$ and $\beta$ as
\begin{equation}
\int_0^{\infty}\frac{e^{\alpha t}-e^{\beta t}}{t} \, dt = 
\ln\left(\frac{\beta}{\alpha}\right).
\label{eqn-2}
\end{equation}
To prove (\ref{eqn-2}), make the substitution $u=e^{-t}$ that turns the 
integral into
\begin{eqnarray}
\int_0^1\frac{u^{-1-\beta}-u^{-1-\alpha}}{\ln u}du & = & 
\int_0^1du\int_{-1-\alpha}^{-1-\beta}u^wdw \nonumber \\
& = & \int_{-1-\alpha}^{-1-\beta}dw\int_0^1u^wdu \nonumber \\
& = & \int_{-1-\alpha}^{-1-\beta}\frac{dw}{w+1} \nonumber \\
& = & \ln\left(\frac{\beta}{\alpha}\right). \nonumber
\end{eqnarray}
\noindent
Now observe that $\vert\frac{2a\cos(rx)}{1+a^2}\vert\leq 1$, therefore it 
is legitimate to expand the logarithmic terms as infinite series using 
$\ln(1+z)=\sum_k(-1)^{k-1}\frac{z^k}k$. The outcome reads
\begin{multline}
\int_0^{\infty}\frac{dx}x\sum_{k\geq1}
\frac{(-1)^{k-1}A^k(\cos^kpx-\cos^kqx)}{k}= \\
\sum_{k\geq1}\frac{(-1)^{k-1}A^k}{2^kk}
\int_0^{\infty}\frac{(e^{ipx}+e^{-ipx})^k
-(e^{iqx}+e^{-iqx})^k}xdx;
\nonumber
\end{multline}
\noindent
where $A=2a/(1+a^2)$. The inner integral is evaluated using some 
 binomial expansions. That is,
\begin{equation}
\int_0^{\infty}\frac{(e^{ipx}+e^{-ipx})^k-(e^{iqx}+e^{-iqx})^k}xdx=
\sum_{r=0}^k\binom{k}{r}\int_0^{\infty}
\frac{e^{(2r-k)ipx}-e^{(2r-k)iqx}}xdx.\label{eqn-3}
\end{equation}
\noindent
It is time to employ equation (\ref{eqn-2}). A closer look at (\ref{eqn-3}) 
shows that care must be exercised. The integrals are sensitive to the 
\it parity \rm of $k$. More precisely, the quantity $2r-k$ vanishes 
if and only if $k$ is even and $r=k/2$, in which case there is a zero 
contribution to summation. Otherwise, the second integral in (\ref{eqn-3}) is 
always equal to $\ln(q/p)$. Therefore,
\begin{equation}
\nonumber
\sum_{r=0}^k\binom{k}{r}
\int_0^{\infty}\left[ e^{(2r-k)ipx}-e^{(2r-k)iqx} \right] 
\frac{dx}{x} =\begin{cases}
2^{k} \ln\left(\frac{q}p\right)\qquad{}\qquad{}\qquad{}\qquad & \text{if $k$ is odd,}\\
\left(2^{k}-\binom{k}{k/2}\right)\ln\left(\frac{q}p\right)
\qquad & \text{if $k$ is even.}\end{cases}
\end{equation}
Combining the results obtained thus far yields 
\begin{eqnarray}
 & & \label{eqn-4} \\
I &= & \int_0^{\infty}\frac{\ln(1+2a\cos(px)+a^2)-\ln(1+2a\cos(qx)+a^2)}{x}dx 
\nonumber \\
& = & 
\int_0^{\infty}\frac{dx}{x}\sum_{k\geq1}\frac{(-1)^{k-1}A^k(\cos^kpx-\cos^kqx)}{k} \nonumber \\
&=& \sum_{k\geq1}\frac{(-1)^{k-1}A^k}{k2^{k}}\sum_{r=0}^k\binom{k}{r}\int_0^{\infty}\frac{e^{(2r-k)ipx}-e^{(2r-k)iqx}}{x}dx \nonumber\\
&=&\ln\left(\frac{q}p\right)
\sum_{\text{$k$ odd}}\frac{(-1)^{k-1}A^k}{k}+
\ln\left(\frac{q}p\right)
\sum_{\text{$k$ even}}\frac{(-1)^{k-1}A^k}{k}\left(1-\frac1{2^k}\binom{k}{k/2}\right) \nonumber \\
&=&\ln\left(\frac{q}p\right)\sum_{k\geq1}\frac{(-1)^{k-1}A^k}k+
\ln\left(\frac{q}p\right)
\sum_{k\geq1}\frac1{2k}\left(\frac{A}{2}\right)^{2k}\binom{2k}k \nonumber \\
&=&\ln\left(\frac{q}p\right)\ln(1+A)+
\frac12\ln\left(\frac{q}p\right)
\sum_{k\geq1}\binom{2k}k\frac1{k}\left(\frac{A^2}{2^2}\right)^{k}.
\nonumber
\end{eqnarray}
\noindent
The last step utilizes the Taylor series
\begin{equation}
\sum_{k\geq1}\binom{2k}k\frac{Q^k}k=
-2\ln\left(\tfrac{1}{2} \left[1 + \sqrt{1-4Q} \right]\right)
\label{eqn-5}
\end{equation}
This follows from the binomial series 
$\sum_{k\geq0}\binom{2k}kR^k=1/\sqrt{1-4R}$ after rearranging in the manner
$$\sum_{k\geq1}\binom{2k}kR^{k-1}=\frac1{R\sqrt{1-4R}}-\frac1{R}=\frac4{\sqrt{1-4R}(1+\sqrt{1-4R})},$$
and then integrating by parts (from $0$ to $Q$)
$$\sum_{k\geq1}\binom{2k}k\frac{Q^k}k=\int_0^Q\frac{4\cdot dR}{\sqrt{1-4R}
(1+\sqrt{1-4R})}=\int_1^{\sqrt{1-4Q}}\frac{-2\cdot du}{1+u}=
-2\ln\left(\tfrac12\left[ 1 + \sqrt{1-4Q} \right]\right).$$
Formula (\ref{eqn-5}) applied to equation (\ref{eqn-4}) leading to 
$$I=\ln\left(\frac{q}p\right)\ln(1+A)-\ln\left(\frac{q}p\right)
\ln\left(\tfrac12 \left[ 1 + \sqrt{1-4Q} \right] \right). $$
\noindent
It remains to replace $Q=A^2/2^2=a^2/(1+a^2)^2$ and use the 
identity
\begin{equation}
1 - 4Q = \frac{(a^{2}-1)^{2}}{(a^{2}+1)^{2}}. 
\nonumber
\end{equation}
\noindent
Observe that the expression for $\sqrt{1-4Q}$  depends on 
whether $|a| > 1$ or not. 
The proof is complete. 

\medskip

\noindent
{\bf Acknowledgments}. Matthew Albano and Erin Beyerstedt were partially 
supported as students by 
$\text{NSF-DMS } 0713836$. The work of the last author was also partially 
supported by the same grant.

\end{document}